\documentclass[12pt]{article}
\usepackage{amsthm,amsfonts,amssymb,amscd}

\begin{document}


\begin{center}
\large \bf Rationally connected Galois covers \\
of Fano-Mori fibre spaces
\end{center}\vspace{0.5cm}

\centerline{A.V.Pukhlikov}\vspace{0.5cm}

\parshape=1
3cm 10cm \noindent {\small \quad\quad\quad \quad\quad\quad\quad
\quad\quad\quad {\bf }\newline For Fano-Mori fibre spaces
$\pi\colon V\to S$, every fibre of which is a primitive Fano
variety with the global canonical threshold at least 1, and which
are stable with respect to birational modifications of the base
and sufficiently twisted over the base, we prove that every
rationally connected Galois cyclic rational cover is pulled back
from the base; in particular, if $S$ does not have such covers
then $V$ does not have rationally connected Galois rational
covers, the Galois group of which is not perfect.

Bibliography: 9 items.}


\vspace{0.3cm}

\section*{Introduction}

{\bf 0.1. Primitive Fano varieties.} It was shown in
\cite{Pukh2020a,Pukh2021a} that primitive divisorially canonical
Fano varieties do not admit rationally connected cyclic Galois
rational covers (which implies immediately that they do not admit
rationally connected rational Galois covers, the Galois group of
which is not perfect, that is, it is not equal to its commutant).
The aim of the present paper is to investigate the problem of
existence of such covers for varieties fibred into primitive Fano
varieties over a positive-dimensional base. In order to state the
main result, recall the necessary definitions.

We work with varieties defined over the field of complex numbers.
Recall that a projective variety $F$ of dimension $M\geqslant 3$
is a {\it primitive Fano variety}, if it is factorial, has at most
terminal singularities and its anticanonical class $(-K_F)$ is
ample and generates the Picard group, $\mathop{\rm Pic} F={\mathbb
Z} K_F$.

A primitive Fano variety $F$ is {\it divisorially canonical}, if
for every effective divisor $D\sim -nK_F$, where $n\geqslant 1$,
the pair $(F,\frac{1}{n}D)$ is canonical: for every exceptional
divisor $E$ over $F$ the inequality
$$
\mathop{\rm ord\,}\nolimits_E D\leqslant n\cdot a(E)
$$
holds (this inequality is opposite to the Noether-Fano
inequality).

In \cite{Pukh2020a,Pukh2021a} the following fact was shown.

{\bf Theorem 0.1.} {\it For a divisorially canonical primitive
Fano variety $F$ and any prime number $p\geqslant 2$ there are no
Galois rational covers $X\stackrel{p:1}{\dashrightarrow} F$, the
Galois group of which is the cyclic group of order $p$, where $X$
is a rationally connected variety. Therefore there are no Galois
rational covers $X\dashrightarrow F$, the Galois group $G$ of
which is not perfect ($G\neq [G,G]$), where $X$ is a rationally
connected variety.}

(Theorem 0.1 is a union of Theorems 1 and 2 in \cite{Pukh2021a}.)
The aim of the present paper is to investigate rationally
connected cyclic Galois rational covers for the total spaces of
fibrations $V\to S$, every fibre of which is a divisorially
canonical primitive Fano variety.\vspace{0.3cm}


{\bf 0.2. Fano-Mori fibre spaces.} A {\it Fano-Mori fibre space}
is a surjective morphism of projective varieties
$$
\pi\colon V\to S,
$$
where $\dim V\geqslant 3 + \dim S$, the base $S$ is non-singular
and rationally connected, and the following conditions are
satisfied:

(FM1) every scheme-theoretic fibre $F_s=\pi^{-1}(s)$, $s\in S$, is
a reduced irreducible factorial Fano variety of dimension $\dim
F=\dim V - \dim S$ with terminal singularities and the Picard
group $\mathop{\rm Pic} F_s\cong {\mathbb Z}$,

(FM2) the variety $V$ itself is factorial and its singularities
are terminal,

(FM3) the equality
$$
\mathop{\rm Pic} V={\mathbb Z} K_V\oplus\pi^* \mathop{\rm Pic}S
$$
holds.

Therefore, Fano-Mori fibre spaces are Mori fibre spaces,
satisfying some additional very good properties.

{\bf Definition 0.1.} A Fano-Mori fibre space $\pi\colon V\to S$
is {\it stable with respect to fibre-wise birational
modifications}, if for every birational morphism $\sigma_S\colon
S^+\to S$, where $S^+$ is a non-singular projective variety, the
morphism
$$
\pi_+\colon V^+=V\mathop{\times}\nolimits_S S^+\to S^+
$$
is a Fano-Mori fibre space.

In the present paper, as in \cite{Pukh2023a}, we study Fano-Mori
fibre spaces, every fibre of which is a divisorially canonical
variety. The proof of the main result (it is stated below in
Subsection 0.3) makes use of slightly weaker properties of fibres.
It is convenient to express them in terms of {\it global
thresholds}.

Recall that for a Fano variety $F$ with the Picard number 1 and
terminal ${\mathbb Q}$-factorial singularities its global
canonical threshold $\mathop{\rm ct}(F)$ is the supremum of
$\lambda\in{\mathbb Q}_+$ such that for every effective divisor
$D\sim -nK_F$ (here $n\in{\mathbb Q}_+$) the pair
$\left(F,\frac{\lambda}{n}D\right)$ is canonical. Therefore, the
divisorial canonicity of a variety $F$ means that the inequality
$\mathop{\rm ct} (F)\geqslant 1$ holds.

If in the definition of the global canonical threshold instead of
``for every effective divisor $D\sim -nK_F$'' we put ``for a
general divisor $D$ in every linear system $\Sigma\subset|-nK_F|$
with no fixed components'', we get the definition of the {\it
mobile canonical threshold} $\mathop{\rm mct}(F)$; obviously,
$\mathop{\rm mct}(F)\geqslant\mathop{\rm ct}(F)$. The inequality
$\mathop{\rm mct}(F)\geqslant 1$ is equivalent to the birational
superrigidity of the Fano variety $F$, see \cite{Ch05c}. If in the
definition of the global canonical threshold the property of the
pair $(F,\frac{\lambda}{n}D)$ to be canonical we replace by the
log canonicity of that pair, we get the definition of the {\it
global log canonical threshold} $\mathop{\rm lct}(F)$; again,
$\mathop{\rm lct}(F)\geqslant\mathop{\rm ct}(F)$.

In \cite{Pukh2023a} the following fact was shown.

{\bf Theorem 0.2.} {\it Assume that a Fano-Mori fibre space
$\pi\colon V\to S$ is stable with respect to fibre-wise birational
modifications and

{\rm (i)} for every point $s\in S$ the fibre $F_s$ satisfies the
inequalitites $\mathop{\rm lct} (F_s)\geqslant 1$ and $\mathop{\rm
mct} (F_s)\geqslant 1$,

{\rm (ii)} (the $K$-condition) every mobile (that is to say, with
no fixed components) linear system on $V$ is a subsystem of a
complete linear system $|-nK_V+\pi^* Y|$, where $Y$ is a
pseudo-effective divisorial class on $S$,

{\rm (iii)} for every family ${\overline{\cal C}}$ of irreducible
curves on $S$, sweeping out a dense subset of the base $S$, and
$\overline{C}\in{\overline{\cal C}}$ no positive multiple of the
class
$$
-(K_V\cdot \pi^{-1}(\overline{C}))-F\in A^{\dim S} V,
$$
where $A^iV$ is the numerical Chow group of classes of cycles of
codimension $i$ on $V$ and $F$ is the class of a fibre of the
projection $\pi$, is represented by an effective cycle on $V$.

Then for every rationally connected fibre space $V'\slash S'$
every birational map $\chi\colon V\dashrightarrow V'$ (if there
are such maps) is fibre-wise, and the fibre space $V\slash S$ is
birationally rigid.}

(This is \cite[Theorem 0.2]{Pukh2023a}.)

In \cite{Pukh15a,Pukh2022a,Pukh2023a} large classes of fibre
spaces satisfying the conditions of Theorem 0.2 were constructed.
One can say that the conditions (ii) and (iii) are typical. The
condition (i) is more restrictive; it is a condition of general
position for $V$ in its family.\vspace{0.3cm}


{\bf 0.3. Galois rational covers.} Let $\pi\colon V\to S$ be a
Fano-Mori fibre space. In this paper we study the rational Galois
covers $\sigma\colon X\dashrightarrow V$ with the cyclic Galois
group of a prime order $p$. In the case when the base $S$ is
positive-dimensional, there is a natural construction of such
covers. Namely, let $\sigma_S\colon Z\dashrightarrow S$ be a
rationally connected rational Galois cyclic cover. Then the fibred
product
$$
V\mathop{\times}\nolimits_S Z
$$
is rationally connected (the fibre of the projection onto the
second factor over a point $s_Z\in Z$ is the fibre of the original
fibre space $\pi\colon V\to S$ over the point $s=\sigma_Z(s_Z)$)
and the projection onto the first factor is a rational cyclic
Galois cover. In this case one can say that the rational Galois
cover is {\it pulled back from the base of the fibre space}. The
main result of this paper is that in the assumptions of Theorem
0.2 the construction described above covers all rationally
connected cyclic Galois rational covers of the variety $V$. (The
notation $V\mathop{\times}\nolimits_S Z$ and the term ``fibred
product'' are meant here in the following ``birational'' sense:
for some open sets $U_Z\subset Z$ and $U_S\subset S$, such that
$\sigma_S|_{U_Z}$ is a cyclic cover $U_Z\to U_S$, we take the
fibred product of the varieties $\pi^{-1}(U_S)\subset V$ and $U_Z$
over $U_S$ and extend it to a projective variety. We use the
notation $V\mathop{\times}\nolimits_S Z$ and the term ``fibred
product'' for simplicity and clarity of exposition.)

{\bf Theorem 0.3.} {\it Assume that a Fano-Mori fibre space
$\pi\colon V\to S$ satisfies all assumptions of Theorem 0.2 and
$\sigma\colon X\dashrightarrow V$ is a cyclic Galois rational
cover of degree $p$, where $X$ is a rationally connected variety.
Then there are a cyclic Galois rational cover $\sigma_S\colon
Z\dashrightarrow S$ of degree $p$, where $Z$ is rationally
connected, and a birational map $\chi\colon X\dashrightarrow
V\mathop{\times}\nolimits_S Z$, such that $\sigma=\mathop{\rm
pr}\nolimits_V\circ\,\chi$, where $\mathop{\rm
pr}\nolimits_V\colon V\mathop{\times}\nolimits_S Z\dashrightarrow
V$ is the projection onto the first factor.}

Setting $\pi_X=\mathop{\rm pr}\nolimits_Z\circ\,\chi$, where
$\mathop{\rm pr}\nolimits_Z$ is the projection of the fibred
product onto the second factor, we get a structure $\pi_X\colon
X\dashrightarrow Z$ of a Fano-Mori fibre space on $X$, where the
diagram
$$
\begin{array}{rcccl}
   & X & \stackrel{\sigma}{\dashrightarrow} & V & \\
\pi_X\!\!\!\!\! & \downarrow &   &   \downarrow & \!\!\!\!\!\pi \\
   & Z & \stackrel{\sigma_S}{\dashrightarrow} & S
\end{array}
$$
is commutative (the projection $\pi_X$, and similarly the other
maps in that diagram can be assumed to be regular), so that the
fibre of general position of the projection $\pi_X$ is birational
to the corresponding fibre of the fibre space $V/S$.

{\bf Corollary 0.1.} {\it Assume that a Fano-Mori fibre space
$\pi\colon V\to S$ satisfies all assumptions of Theorem 0.2 and
the variety $S$ has no rationally connected Galois rational covers
of the prime order $p$. Then the variety $V$ also has no
rationally connected Galois rational covers of order $p$.}

{\bf Corollary 0.2.} {\it Assume that a Fano-Mori fibre space
$\pi\colon V\to S$ satisfies all assumptions of Theorem 0.2 and
the variety $S$ has no rationally connected Galois rational
covers, the Galois group $G$ of which is not perfect. Then the
variety $V$ also has no such Galois covers.}

Corollary 0.2 allows one to iteratively construct families of
rationally connected varieties that have no rationally connected
rational Galois covers, the Galois group of which is not perfect,
see below.\vspace{0.3cm}


{\bf 0.4. Explicit examples.} Let us consider the most visual
class of fibrations into Fano hypersurfaces of index 1. It is
known \cite[Theorem 1.4]{Pukh15a}, that for $M\geqslant 10$ in the
projective space ${\cal F}={\mathbb P}(H^0({\mathbb P}^M,{\cal
O}_{{\mathbb P}^M}(M)))$, parameterizing hypersurfaces of degree
$M$ in ${\mathbb P}^M$, there is an open subset ${\cal F}_{\rm
reg}$, such that for every point of that subset the corresponding
hypersurface is irreducible, reduced factorial and divisorially
canonical primitive Fano variety with terminal singularities, and
moreover, the codimension of the complement ${\cal F}\setminus
{\cal F}_{\rm reg}$ is at least $\frac12(M-7)(M-6)-5$. Now let $S$
be a non-singular rationally connected projective variety of
dimension $\leqslant\frac12 (M-7)(M-6)-6$ and ${\cal L}$ a locally
free sheaf of rank $M+1$ on $S$, generated by global sections, so
that ${\cal X}={\mathbb P}({\cal L})={\bf Proj}\,
\mathop{\oplus}\limits_{i=0}^{\infty}{\cal L}^{\otimes i}$ is a
locally trivial ${\mathbb P}^M$-fibration over $S$. Set $\pi_{\cal
X}\colon {\cal X}\to S$ to be the projection and let $L\in
\mathop{\rm Pic} {\cal X}$ be the class of the sheaf ${\cal
O}_{{\mathbb P}({\cal L})}(1)$. Now we construct a Fano-Mori fibre
space $\pi\colon V\to S$ in a natural way: $V\in |ML+\pi^*_{\cal
X} R|$, where $R\in \mathop{\rm Pic} S$ is some class, is a
sufficiently general divisor, so that for every point $s\in S$ we
have
$$
\pi^{-1}(s)\in {\cal F}_{\rm reg}.
$$
Now \cite[Theorem 1.5]{Pukh15a} states that if
$K_S+\left(1-\frac{1}{M}\right)R$ is a pseudo-effective ${\mathbb
Q}$-divisorial class, then the fibre space $\pi\colon V\to S$
satisfies all assumptions of Theorem 0.2 (the stability with
respect to fibre-wise birational modifications follows from the
explicit description of singularities of the hypersurfaces $F\in
{\cal F}_{\rm reg}$, see \cite[Subsection 0.3]{Pukh2023a} in a
more general context of fibrations into complete intersections).
Therefore for the fibre space $V/S$ the claim of Theorem 0.3 is
true.

{\bf Example 0.1.} Take $m\leqslant\frac12 (M-7)(M-6)-6$ and let
$S \subset {\mathbb P}^{m+1}$ be a non-singular divisorially
canonical hypersurface of degree $m+1$, ${\cal X}={\mathbb
P}^M\times S$, $L$ the pull back of the class of a hyperplane in
${\mathbb P}^M$ on ${\cal X}$ and $V\in|ML+\pi^*_{\cal X} R|$,
where $R=lH_S$ (here $H_S$ is the class of a hyperplane section of
the hypersurface $S$). It is easy to check that the condition for
the class $K_S+\left(1-\frac{1}{M}\right)R$ to be pseudo-effective
is equivalent to the inequality $l\geqslant 2$. Therefore, if
$(x_0:\dots :x_M)$ and $(y_0:\dots :y_{m+1})$ are homogeneous
coordinates on ${\mathbb P}^M$ and ${\mathbb P}^{m+1}$,
respectively, then the complete intersection $V$ in ${\mathbb
P}^M\times{\mathbb P}^{m+1}$, given by a pair of sufficiently
general polynomials $f_S(y_0,\dots,y_{m+1})$ of degree $m+1$ and
$f_V(x_0,\dots, x_M;y_0,\dots, y_{m+1})$ of bi-degree $(M,l)$,
where $l\geqslant 2$, does not admit rationally connected Galois
rational covers, the Galois group of which is not perfect (in
particular, rational double covers).\vspace{0.3cm}


{\bf 0.5. The structure of the paper and general remarks.} The
purpose of this paper is to prove Theorem 0.3. We do it in three
steps. First (\S 1), we construct a cyclic cover
$$
\sigma_U\colon U_X\to U
$$
of non-singular quasi-projective varieties, where $U_X$ is
birational to $X$ and $U$ is a Zariski open subset of some
birational modification $\widetilde{V}$ of the variety $V$, which
is equivalent to the original Galois rational cover
$X\dashrightarrow V$ (that is, the corresponding extensions of the
fields of rational functions are the same). Although $U_X$ and $U$
are quasi-projective, they contain families of projective rational
curves. The cover $\sigma_U$ is branched over a non-singular (but
possibly reducible and in that case non-connected) hypersurface
$W\subset U$. The image on $V$ of the union of the components of
the divisor $W$, which are divisorial on $V$ and cover the base
$S$, has the class $-nK_V$ up to some divisor that is pulled back
from the base.

After that (\S 2), we exclude the option $n\geqslant 2$. This is
done by the same technique that was used in \cite{Pukh2023a} and
the previous papers \cite{Pukh15a,Pukh2022a} to prove the
birational rigidity (see Theorem 0.2). Therefore, one of the two
cases, $n=0$ or $n=1$ takes place.

Finally, in \S 3 we exclude the case $n=1$, and in the case $n=0$
prove that the cover $\sigma\colon X\dashrightarrow V$ is pulled
back from the base (this follows from the fact that the
``essential'' part of the divisor $W$ is pulled back from the
base), which completes the proof of Theorem 0.3.

The paper makes use of the technique of the method of maximal
singularities (see, for instance, \cite[\S 2]{Pukh2023a}), some
well known facts about free families of rational curves on
rationally connected varieties \cite[Chapter II]{Kol96} and the
inversion of adjunction, based on the Shokurov-Koll\'{a}r
connectedness principle \cite{Kol93,Sh92}.

The author is grateful to the members of Divisions of Algebraic
Geometry and Algebra at Steklov Institute of Mathematics for the
interest to his work, and to the colleagues in Algebraic Geometry
research group at the University of Liverpool for general support.


\section{The cyclic cover}

In this section we construct a cyclic cover $\sigma_U\colon U_X\to
U$ of non-singular quasi-projective varieties, which is equivalent
(in the sense of the extension of the fields of rational
functions) to the original Galois rational cover $\sigma\colon
X\dashrightarrow V$. The branch divisor $W\subset U$ accumulates
the key information about the cover $\sigma_U$.\vspace{0.3cm}

{\bf 1.1. The mobile family of curves.} Let us fix a rational
Galois cyclic cover $\sigma\colon X\dashrightarrow V$ of degree
$p$, where $X$ is a rationally connected variety. Without loss of
generality we assume that $X$ is a non-singular projective variety
and $\sigma$ is a morphism. Let us consider a family ${\cal C}_X$
of non-singular irreducible rational curves on $X$, satisfying the
following properties:

(C1) the curves $C_X\in{\cal C}_X$ sweep out a dense subset of the
variety $X$,

(C2) for a general pair of points on $X$ there is a curve
$C_X\in{\cal C}_X$ containing them,

(C3) for every subvariety $Y\subset X$ of codimension $\geqslant
2$ the set of curves $C_X\in {\cal C}_X$, meeting $Y$, forms a
proper closed subfamily of the family ${\cal C}_X$ (in other
words, a curve $C_X$ of general position does not meet the subset
$Y$),

(C4) for every prime divisor $\Delta\subset X$, such that either
$\mathop{\rm codim}(\sigma(\Delta)\subset V)\geqslant 2$, or for a
point $a\in\Delta$ of general position we have $\sigma_*\colon
T_aX\to T_{\sigma(a)}V$ is not an isomorphism, a general curve
$C_X\in{\cal C}_X$ either does not meet $\Delta$, or intersects
$\Delta$ transversally at a set of points of general position,

(C5) for every curve $C_X\in{\cal C}_X$ the morphism
$\sigma|_{C_X}\colon C_X\to\sigma(C_X)$ is birational.

The existence of such family is easy to show by means of
deformation arguments \cite[Chapter II]{Kol96}, see
\cite{Pukh2020a} for the details. Since we will only need a curve
of general position in that family, we may remove from ${\cal
C}_X$ some proper closed subfamilies when we need it (for
instance, the subfamily of curves meeting a given closed subset of
codimension $\geqslant 2$) without special comments and keeping
the notation ${\cal C}_X$ for the smaller family.

A family of irreducible projective rational curves on a
quasi-projective variety, satisfying the properties (C1-3), is
said to be {\it free.}

Set ${\cal C}_V=\sigma_*{\cal C}_X$. It is a family of irreducible
rational curves on $V$. By (C2) these curves are not contracted by
the projection $\pi$. Set ${\cal C}_S=\pi_*{\cal C}_V$. This is a
family of irreducible rational curves on $S$.

By deformation arguments, based on \cite[Sections II.1 and
II.3]{Kol96}, see \cite[Section 3]{Pukh2020a} for the details, we
may assume that for a general (and thus for every) curve
$C_X\in{\cal C}_X$
$$
\pi\circ\sigma|_{C_X}\colon C_X\to C_S=\pi\sigma(C_X)
$$
is a birational morphism. However, the family ${\cal C}_S$,
generally speaking, is not free: of course, the properties (C1,2)
are satisfied, but the property (C3) may not be true; for
instance, if some prime divisor on $X$ that has a non-empty
intersection with a general curve $C_X$ is contracted by the map
$\pi\sigma$ to a subvariety of codimension $\geqslant 2$ on
$S$.\vspace{0.3cm}


{\bf 1.2. Resolution of a family of curves.} Let
$$
\lambda_S\colon S^+\to S
$$
be a resolution of the family ${\cal C}_S$ in the sense of
\cite[Proposition 1]{Pukh2020a}: the variety $S^+$ is non-singular
and projective and the strict transform ${\cal C}^+_S$ of the
family ${\cal C}_S$ is a free family of irreducible rational
curves on $S^+$. (For a proof that a resolution does exist, see
\cite[Section 4]{Pukh2020a}.)

Set $V^+=V\times_SS^+$ and denote by the symbols $\lambda$ and
$\pi_+$ the projections of $V^+$ onto $V$ and $S^+$, respectively.
By the assumptions about the fibre space $V\slash S$ the
projection $\pi_+\colon V^+\to S^+$ is a Fano-Mori fibre space
over $S^+$, see \cite[Subsection 2.2]{Pukh2023a}. Let ${\cal
C}^+_V$ be the strict transform of the family ${\cal C}_V$ on
$V^+$. Obviously, $(\pi_+)_*{\cal C}^+_V={\cal C}^+_S$.

The family of curves ${\cal C}^+_V$ may be not free. Let
$\varphi\colon\widetilde{V}\to V^+$ be its resolution (a sequence
of blow ups with non-singular centers), so that $\widetilde{V}$ is
a non-singular projective variety and the strict transform
$\widetilde{\cal C}_V$ of the family ${\cal C}^+_V$ on
$\widetilde{V}$ is a free family of irreducible rational curves on
$\widetilde{V}$.

{\bf Proposition 1.1.} {\it For every irreducible subvariety
$Y\subset S^+$ of codimension $\geqslant 2$ the set of curves
$\widetilde{C}_V\in\widetilde{\cal C}_V$, such that
$(\pi_+\circ\varphi)\widetilde{C}_V\cap Y\neq\emptyset$ is a
proper closed subfamily of the family} $\widetilde{\cal C}_V$.

{\bf Proof.} By construction,
$$
(\pi_+)_*\varphi_*\widetilde{\cal C}_V={\cal C}^+_S,
$$
and the family ${\cal C}^+_S$ is free. Q.E.D. for the
proposition.\vspace{0.3cm}


{\bf 1.3. The cyclic cover.} Now everything is ready to construct
the cyclic cover, equivalent to the Galois rational cover
$\sigma\colon X\to V$.

{\bf Proposition 1.2.} {\it There exist a non-singular
quasi-projective variety $U_X$, a birational map $\varphi_X\colon
U_X\dashrightarrow X$ and a Zariski open subset
$U\subset\widetilde{V}$, such that:

{\rm (i)} the rational map
$$
\varphi^{-1}\circ\lambda^{-1}\circ\sigma\circ\varphi_X\colon
U_X\dashrightarrow\widetilde{V}
$$
extends to a morphism $\sigma_U\colon U_X\to\widetilde{V}$, the
image of which is $U$,

{\rm (ii)} $\mathop{\rm codim}((\widetilde{V}\backslash
U)\subset\widetilde{V})\geqslant 2$,

{\rm (iii)} the morphism $\sigma_U\colon U_X\to U$ is a cyclic
cover of degree $p$, branched over a non-singular hypersurface}
$W\subset U$.

{\bf Proof} repeats the arguments in \cite[Section 5]{Pukh2020a}
word for word, however, certain notations have to be changed: for
instance, the symbol $\widetilde{V}$ in the present paper has the
same meaning as the symbol $V^+$ in \cite{Pukh2020a}, whereas the
latter symbol in the present paper is used for another purpose. In
order to avoid any confusion, we will trace these changes in the
notations.

Starting from the cyclic field extension ${\mathbb C}(V)={\mathbb
C}(\widetilde{V})\subset{\mathbb C}(X)$, we construct an
irreducible hypersurface $X_0$ in the direct product
$\widetilde{V}\times{\mathbb P}^1_{(x_0 : x_1)}$, given by the
equation
$$
a_1x_1^p-a_0x^p_0=0,
$$
where $(a_1x_1^p-a_0x^p_0)$ is a section of the sheaf
$$
\mathop{\rm pr}\nolimits^*_V{\cal
O}_{\widetilde{V}}(R)\otimes\mathop{\rm pr}\nolimits^*_{\mathbb
P}{\cal O}_{{\mathbb P}^1}(p),
$$
the symbols $\mathop{\rm pr}_V$ and $\mathop{\rm pr}_{\mathbb P}$
denote the projections onto the first and second factor,
respectively, $R$ is some effective divisor on $\widetilde{V}$ and
the sections $a_0,a_1\in{\cal O}_{\widetilde{V}}(R)$ do not vanish
simultaneously on any prime divisor on $\widetilde{V}$, and
moreover, $X_0$ is birational to the original variety $X$: there
is a commutative diagram of maps
$$
\begin{array}{ccc}
X_0 & \dashrightarrow & X\\
\downarrow & & \downarrow\\
\widetilde{V} & \stackrel{\lambda\circ\varphi}{\longrightarrow} & V,\\
\end{array}
$$
where the upper horizontal arrow is a birational map and the left
vertical arrow is induced by the projection $\mathop{\rm pr}_V$.

Let $|a_0=0|$ and $|a_1=0|$ be the supports of the divisors of
zeros of the sections $a_0$, $a_1$, that is,
$$
|a_i=0|=\{q\in\widetilde{V}\,|\,a_i=0\}.
$$
The open subset $U\subset\widetilde{V}$ is obtained by removing
from $\widetilde{V}$ the closed subset
$$
|a_0=0|\cap|a_1=0|
$$
and all singularities of the hypersurfaces $|a_0=0|$ and
$|a_1=0|$. By construction, the hypersurfaces $|a_0|_U=0|$ and
$|a_1|_U=0|$ on $U$ are disjoint and non-singular (but possibly
not connected). Now we argue in the word for word the same way as
in \cite[Section 5]{Pukh2020a} (see also \cite[Section
6]{Pukh2021a}) and complete the proof of the proposition.

Let $\Theta\subset U$ be a prime divisor. At least one of the
integers $\mathop{\rm ord}_{\Theta}a_0$, $\mathop{\rm
ord}_{\Theta}a_1$ is zero. Set
$$
\mu(\Theta)=\mathop{\rm max}\{\mathop{\rm
ord}\nolimits_{\Theta}a_0,\,\mathop{\rm
ord}\nolimits_{\Theta}a_1\}.
$$

{\bf Proposition 1.3.} {\it The hypersurface  $W$ contains
$\Theta$ if and only if} $p\not |\mu(\Theta)$.

{\bf Proof}: this is \cite[Proposition 1]{Pukh2021a}, a proof is
given in \cite[Section 6]{Pukh2021a}.

Let $\overline{W}$ be the closure of $W$ in $\widetilde{V}$, so
that $\overline{W}\cap U=W$. Write down
\begin{equation}\label{14.09.23.1}
\lambda_*\varphi_*\overline{W}|_{F_s}\sim-nK_{F_s},
\end{equation}
where $F_s=\pi^{-1}(s)$ is a fibre of general position,
$K_{F_s}=K_V|_{F_s}$ its canonical class and $n\in{\mathbb Z}_+$.
There are the following options:

(W0) $n=0$,

(W1) $n=1$,

(W2)  $n\geqslant 2$.

We will show that (W1) and (W2) do not realize and (W0) means that
$\sigma$ is pulled back from the base.


\section{Exclusion of the case $n\geqslant 2$}

In this section we show that the case (W2), that is, when
$n\geqslant 2$, is impossible, so that
$n\in\{0,1\}$.\vspace{0.3cm}

{\bf 2.1. Divisorial classes on $S^+$, $V^+$ and $\widetilde{V}$.}
Let ${\cal T}$ be the set of all prime $\lambda_S$-exceptional
divisors on $S^+$ and ${\cal E}$ the set of all prime
$\varphi$-exceptional divisors on $\widetilde{V}$. By the
properties of the Fano-Mori fibre space $\pi\colon V\to S$ and the
construction of the variety $V^+$ we know that
$\lambda$-exceptional divisors on $V^+$ are precisely all divisors
of the form $\pi^*_+T$, where $T\in{\cal T}$, and for the
discrepancies we have the equality
$$
a(\pi^*_+T,V)=a(T,S).
$$
Denote the last discrepancy by the symbol $a_T$. Similarly, we
denote the discrepancy of $E\in{\cal E}$ with respect to $V^+$ by
the symbol $a_E$. Now for the canonical classes of the varieties
$S^+$, $V^+$ and $\widetilde{V}$ we get:
$$
K^+_S=\lambda^*_SK_S+\sum_{T\in{\cal T}}a_TT,\quad
K^+_V=\lambda^*K_V+\pi^*_+\left(\sum_{T\in{\cal T}}a_TT\right)
$$
and
$$
K_{\widetilde{V}}=\varphi^*K^+_V+\sum_{E\in{\cal E}}a_EE,
$$
respectively; all discrepancies are positive integers. (The
notations $K^+_S$ and $K^+_V$ for the canonical classes of $S^+$
and $V^+$ are non-standard but convenient.) Finally, by the symbol
$K$ denote the canonical class of the variety $U_X$. By
construction,
$$
K=\sigma^*_U\left(K_{\widetilde{V}}|_U+\left(1-\frac{1}{p}\right)W\right).
$$

Let us write down $\overline{W}=W_{\rm div}+W_{\rm exc}$, where
all components in $W_{\rm div}$ are divisorial on $V$, and $W_{\rm
exc}$ consists of $\lambda\circ\varphi$-exceptional divisors.
Obviously,
$$
W_{\rm exc}=\sum_{E\in{\cal E}}c_EE+\sum_{T\in{\cal
T}}c_T(\pi^{-1}_+(T))^{\sim},
$$
where $c_E,c_T\in\{0,1\}$, and the upper index $\sim$ means the
strict transform on $\widetilde{V}$. Furthermore, set
$$
\lambda_*\varphi_*W_{\rm div}=D+\Delta,
$$
where $D$ is an effective divisor, each component of which covers
the base $S$, and $\Delta$ is an effective divisor, pulled back
from $S$. Taking into account (\ref{14.09.23.1}), we can write
$$
D\sim-nK_V+\pi^*Y,
$$
where $Y$ is a (pseudoeffective) divisorial class on $S$. The
strict transforms of the effective divisors $D$ and $\Delta$ on
$V^+$ and $\widetilde{V}$ we denote by the symbols $D^+$ and
$\Delta^+$, and $\widetilde{D}$, $\widetilde{\Delta}$,
respectively. Therefore, $W_{\rm div}=\widetilde{D} +
\widetilde{\Delta}$ and $\varphi_*W_{\rm div}=D^++\Delta^+$. We
get
$$
D^+\sim -nK_V+\pi^*_+\left(Y-\sum_{T\in{\cal T}}b_TT\right)
$$
and
$$
\widetilde{D}\sim -nK_V+\varphi^*\pi^*_+\left(Y-\sum_{T\in{\cal
T}}b_TT\right)-\sum_{E\in{\cal E}}b_EE,
$$
where $b_T$, $b_E\geqslant 0$ and we write $K_V$ instead of
$\lambda^*K_V$ and $\varphi^*\lambda^*K_V$ and $Y$ instead of
$\lambda^*_SY$ for simplicity. Combining the formulas for $D^+$
and $\widetilde{D}$ with the formulas for the canonical classes,
given above, we see that the class of the divisor $\widetilde{D}$
is
$$
-nK_{\widetilde{V}}+\varphi^*\pi^*_+\left(Y+\sum_{T\in{\cal
T}}(na_T-b_T)T\right)+\sum_{E\in{\cal E}}(na_E-b_E)E.
$$


{\bf 2.2. The curves in the family ${\widetilde{\cal C}}_V$.}
Assume now that $n\geqslant 2$ (the case (W2)). Since the family
of curves $\widetilde{\cal C}_V$ is by construction free, the
claim (ii) of Proposition 1.2 gives us that the curves
$\widetilde{C}\in\widetilde{\cal C}_V$ that are entirely contained
in $U$, form a free family ${\cal C}_U$ of irreducible projective
rational curves on $U$ (and it is an open subfamily of the family
$\widetilde{\cal C}_V$). By the property (C5) for every curve
$C\in{\cal C}_U$ we get:
$$
\sigma^{-1}_U(C)=C_1\cup C_2\cup\dots\cup C_p
$$
is a union of $p$ distinct rational curves, permuted by the
elements of the Galois group, and moreover, $\sigma_U|_{C_i}\colon
C_i\to C$ is birational for all $i=1,\dots,p$. Each curve $C_i$
varies in a family of irreducible projective rational curves,
sweeping out a dense subset in $U_X$ (see \cite[Section
3]{Pukh2020a}), so that $(C_i\cdot K)<0$ (recall that $K$ is the
canonical class of $U_X$).

Since the curves $C\in{\cal C}_U$ sweep out a dense subset of
$\widetilde{V}$, they have a non-negative intersection with every
effective divisor on $\widetilde{V}$. On the other hand, the
inequality $(C_i\cdot K)<0$ for $i=1,\dots,p$ implies that
$$
\left(C\cdot\left(K_{\widetilde{V}}+\left(1-\frac{1}{p}\right)W\right)\right)<0.
$$
Let us add to the left hand side of that inequality the
non-positive (as $n\geqslant 2$) number
$$
\left[\left(1-\frac{1}{p}\right)n-1\right](C\cdot
K_{\widetilde{V}})
$$
and obtain the inequality
$$
n(C\cdot K_{\widetilde{V}})+(C\cdot W)<0,
$$
which only becomes stronger if we replace $W$ by $\widetilde{D}$
(the effective divisor $\widetilde{D}$ is a part of the effective
divisor $W$). Therefore,
\begin{equation}\label{16.09.23.1}
\left(C\cdot \varphi^*\pi_+^*\left(Y+\sum_{T\in{\cal
T}}(na_T-b_T)T\right)\right)+\left( C\cdot\sum_{E\in{\cal
E}}(na_E-b_E)E\right)<0.
\end{equation}

The following fact is of key importance for our arguments.

{\bf Proposition 2.1.} {\it If for $E\in{\cal E}$ the inequality
$(C\cdot E)>0$ holds, then} $b_E\leqslant na_E$.

{\bf Proof.} Since $(C\cdot E)>0$, the curve $(\pi_+\circ
\varphi)(C)\in{\cal C}^+_S$ meets the set $\pi_+\circ
\varphi(E)\subset S^+$. Now it follows from Proposition 1.1 that
$\pi_+\circ \varphi(E)$ is either the whole variety $S^+$, or a
prime divisor on $S^+$. Since each component of the divisor $D$
covers the base $S$, each component of the divisor $D^+$ (the
strict transform of the divisor $D$ on $V^+$) covers $S^+$, so
that $D^+$ does not contain entirely the divisor
$\pi^{-1}_+(\pi_+\circ\varphi(E))$ in the second case. In any
case, for a point of general position $s_+\in\pi_+\circ\varphi(E)$
the corresponding fibre $F=\pi^{-1}_+(s_+)$ is not contained in
$D^+$. Assume now that $b_E>na_E$. Then the pair
$(V^+,\frac{1}{n}D^+)$ is not canonical along $\varphi(E)$, so
that the pair $(F^+,\frac{1}{n}D_F^+)$, where $D^+_F=D^+|_F$, is
not canonical (in the second case even not log canonical by the
inversion of adjunction). This contradicts the condition (i) of
Theorem 0.2. Q.E.D. for the proposition.

We conclude that the second part in the left hand side of the
inequality (\ref{16.09.23.1}) is non-negative, so that for a
general curve $C^+_S\in {\cal C}^+_S$ (recall that ${\cal
C}^+_S=(\pi_+\circ\varphi)_* {\cal C}_U$) the inequality
$$
\left( C^+_S \cdot\left(Y+\sum_{T\in{\cal
T}}(na_T-b_T)T\right)\right)<0
$$
holds. Note that since the class $Y$ is pseudoeffective, we have
the inequality $(C^+_S\cdot Y)\geqslant 0$, so that ${\cal
T}\neq\emptyset$ and for some $T\in{\cal T}$ we have $(C^+_S\cdot
T)>0$, and therefore
$$
\left( C^+_S \cdot\left(Y-\sum_{T\in{\cal T}}b_T
T\right)\right)<-n.
$$
Now we argue word for word as in \cite[Subsection 2.4]{Pukh2023a}:
consider the algebraic cycle of the scheme-theoretic intersection
$$
(D^+\circ \pi^{-1}_+(C^+_S))=\left(\left(\lambda^* D-
\pi^*_+\left(\sum_{T\in{\cal T}}b_T T\right)\right)\circ
\pi^{-1}_+(C^+_S)\right).
$$
This cycle is effective and its numerical class is
$$
n(\lambda^*(-K_V)\cdot \pi^{-1}_+(C^+_S))+\left(C^+_S\cdot
\left(\lambda^*_S Y- \sum_{T\in{\cal T}}b_T T\right)\right)F,
$$
where $F$ is the class of a fibre of the projection $\pi_+$, so
that the class of the effective cycle $\lambda_*(D^+\circ
\pi^{-1}_+(C^+_S))$ in the numerical Chow group of the variety $V$
is
$$
-n(K_V\cdot \pi^{-1}(C_S))+bF,
$$
where $b<-n$ (and the symbol $F$ now stands for the class of a
fibre of the projection $\pi$). This contradicts the condition
(iii) of Theorem 0.2 and completes the exclusion of the case
$n\geqslant 2$.


\section{Reduction to the base of the fibre space}

In this section we exclude the case $n=1$ and prove that if $n=0$,
then the Galois rational cover is pulled back from the base, which
completes the proof of Theorem 0.3.\vspace{0.3cm}

{\bf 3.1. Restriction onto a fibre of general position.} Set
$F=\pi^{-1}_+(s_+)$, where $s_+\in S^+$ is a point of general
position, so that $s=\lambda_S(s_+)\in S$ is a point of general
position, too. Let $\widetilde{F}$ be the strict transform of the
fibre $F$ on $\widetilde{V}$. Set
$$
U_F=F\backslash\varphi\left[\left(\left(\bigcup_{E\in{\cal
E}}E\right)\cup(\widetilde{V}\backslash
U)\right)\cap\widetilde{F}\right].
$$
This is an open subset in $F$, and moreover,
\begin{equation}\label{30.09.23.1}
\mathop{\rm codim}((F\backslash U_F)\subset F)\geqslant 2.
\end{equation}
Obviously, $U_F$ is isomorphic to the open subset
$$
\widetilde{F}\backslash\left[\left(\left(\bigcup_{E\in{\cal
E}}E\right)\cup(\widetilde{V}\backslash
U)\right)\cap\widetilde{F}\right],
$$
which by construction is contained in $U$. We identify these two
quasi-projective varieties, so that it makes sense to write $W\cap
U_F$. Obviously,
$$
W_F=\overline{W\cap U_F}=D\cap F\sim-nK_F
$$
(we use the notations of Subsection 2.1 and identify $F$ with the
fibre $\pi^{-1}(s)$). In the notations of Subsection 1.3 set
$$
R_F=\overline{R\cap U_F}
$$
and let $a_{F,i}\in{\cal O}_F(R_F)$ be the sections $a_i|_{U_F}$,
extended to the whole variety $F$ (this makes sense as the
inequality (\ref{30.09.23.1}) holds and the variety $F$ is
factorial), where $i=0,1$. We have $R_F\sim-NK_F$ for some
$N\geqslant 1$. Recall that $\widetilde{D}$ is the strict
transform of the divisor $D$ on $\widetilde{V}$. Obviously,
$$
\widetilde{D}\cap U_F=W\cap U_F.
$$
Now everything is ready for investigating the cases $n=1$ and
$n=0$. Our main instrument is Proposition 1.3, which makes it
possible to compare the divisors $W_F$ and $R_F$.\vspace{0.3cm}


{\bf 3.2. Exclusion of the case $n=1$.} Assume that $n=1$. In that
case $D$ is a prime divisor (an irreducible hypersurface, covering
$S$) and $W_F\sim-K_F$ is a prime divisor on $F$. Assume for
convenience that $\mathop{\rm ord}_{\widetilde{D}}a_0=0$. By
Proposition 1.3 we have $p\not|\mathop{\rm
ord}_{\widetilde{D}}a_1$. Therefore, $\mathop{\rm
ord}_{W_F}a_{F,0}=0$ and $p\not|\mathop{\rm ord}_{W_F}a_{F,1}$. By
the same Proposition 1.3, for every prime divisor $\Theta\neq W_F$
on $F$ we have
$$
p\,|\,\mathop{\rm ord}\nolimits_{\Theta}a_{F,i},\quad i=0,1.
$$
Therefore, the prime integer $p$ divides the order of the section
$a_{F,0}$ along every prime divisor $\Theta$ on $F$, which implies
that $p\,|\,N$.

However, $W_F\sim-K_F$, whence in a similar way, considering the
section $a_{F,1}$, we get that $p\not|\, N$. This contradiction
excludes the case $n=1$.\vspace{0.3cm}


{\bf 3.3. The cyclic cover of the base.} Assume now that $n=0$. In
that case $D=0$ and $W_F=0$, so that by Proposition 1.3 we get:
$p\,|\mathop{\rm ord}_{\Theta}a_{F,i}$ for every prime divisor
$\Theta$ on $F$ and $i=0,1$. In the notations of Subsection 1.3 we
see that
$$
\mathop{\rm pr}\nolimits^{-1}_V(U_F)\cap X_0\subset
U_F\times{\mathbb P}^1_{(x_0:x_1)}
$$
is given by the equation
\begin{equation}\label{18.09.23.1}
a_{F,1}x^p_1-a_{F,0}x^p_0=0,
\end{equation}
where $a_{F,i}\in{\cal O}_F(-NK_F)$, where $p\,|\,N$ and we can
``extract the root'' from the sections $a_{F,i}$: there are
sections
$$
e_{F,0},e_{F,1}\in H^0(F,{\cal O}_F((-N\slash p)K_F)),
$$
such that $a_{F,i}=e^p_{F,i}$, $i=0,1$. We conclude that the left
hand side of the equation (\ref{18.09.23.1}) is reducible:
$$
a_{F,1}x^p_1-a_{F,0}x^p_0=\prod^p_{i=1}(e_{F,1}x_1-\zeta^ie_{F,0}x_0),
$$
where $\zeta=\mathop{\rm exp}(2\pi i\slash p)$, so that
$$
\sigma^{-1}(F_s)=F_{s,1}\cup F_{s,2}\cup\dots\cup F_{s,p}
$$
is a union of $p$ irreducible components, the map
$\sigma|_{F_{s,i}}\colon F_{s,i}\to F_s$ is birational and the
cyclic Galois group of the extension ${\mathbb
C}(V)\subset{\mathbb C}(X)$ acts transitively on the set of these
components. This is true for the fibres $F_s$ over the points
$s\in U_S$ of some Zariski open subset of the base $U_S\subset S$.

Let $U_Z$ be a non-singular quasi-projective variety,
parameterizing the irreducible components $F_{s,i}$ for $s\in
U_S$, and $Z$ a non-singular projective variety, containing $U_Z$
as a Zariski open subset. Mapping an irreducible component
$F_{s,i}$ to the fibre $F_s$, we obtain a $p$-sheeted Galois
cyclic cover $U_Z\to U_S$, which extends to a Galois rational
cyclic cover
$$
\sigma_S\colon Z\dashrightarrow S
$$
(we can and will assume that $\sigma_S$ is a regular map).
Moreover, there is the obvious rational dominant map
$$
\pi_X\colon X\dashrightarrow Z,
$$
the fibre $\pi^{-1}_X(s_Z)$, $s_Z\in Z$, of which is birational to
the fibre $\pi^{-1}(\sigma_S(s_Z))=F_{\sigma_S(s_X)}$ and which
generates the commutative diagram
$$
\begin{array}{rcccl}
  & X & \stackrel{\sigma}{\longrightarrow} & V &  \\
\pi_X \!\!\! & \downarrow & & \downarrow & \!\!\!\pi\\ &
Z & \stackrel{\sigma_Z}{\longrightarrow} & S. &\\
\end{array}
$$
(As the image of the variety $X$, the variety $Z$ is rationally
connected.) In other words, the rationally connected variety $X$
is birational to the fibred product $V\times_SZ$, and with respect
to that (birational) identification the Galois cyclic rational
cover $\sigma$ is the projection onto the first factor, whereas
$X$ itself has a structure of a Fano-Mori fibre space, defined by
the projection onto the second factor.

Q.E.D. for Theorem 0.3.


\begin{flushleft}
Department of Mathematical Sciences,\\
The University of Liverpool
\end{flushleft}

\noindent{\it pukh@liverpool.ac.uk}

\end{document}